\documentclass[12pt]{article}
\usepackage{epsfig}
\usepackage{latexsym}
\usepackage{amssymb}
\usepackage{amsfonts}
\usepackage{amsmath}
\usepackage{theorem}
\usepackage{graphics}
\newcommand{\bnull}{{\boldsymbol 0}}

\newcommand{\ba}{{\bf a}}

\newcommand{\bp}{{\bf p}}
\newcommand{\bb}{{\bf b}}
\newcommand{\bz}{{\bf z}}
\newcommand{\bx}{{\bf x}}

\newcommand{\by}{{\bf y}}

\newcommand{\bs}{{\bf s}}

\newcommand{\eot}{{\begin{flushright} \vspace{-0.7cm} {$ \blacksquare $} \end{flushright}}}

\newcommand{\lra}{\longrightarrow}
\newcommand{\ra}{\rightarrow}

\newcommand{\beq}{\begin{equation}}
\newcommand{\eeq}{  \end{equation}}

\newcommand{\supp}{\operatorname{supp}}

\newcommand{\R}{\mathbb{R}}

\newcommand{\N}{\mathbb{N}}

\newcommand{\Proof}{{\sc Proof:} \\}

\newcommand{\cC}{{\cal C}}

\newcommand{\cK}{{\cal K}}

\newcommand{\cN}{{\cal N}}
\newcommand{\cO}{{\cal O}}

\newcommand{\cS}{{\cal S}}

\newcommand{\cX}{{\cal X}}

\newtheorem{theorem}{Theorem}

\newtheorem{lemma}[theorem]{Lemma}

\begin{document}
\bibliographystyle{plain}
\parindent = 0.0cm
\setlength{\parskip}{1ex plus 0.5ex minus 0.5ex} \sloppy \frenchspacing
\title{Continuity of scalar-fields characterized by smooth\\ paths fulfilling $\|\bs(t)\|\|\bs'(t)\| < +\infty$.}

\author{
Sigurdur F.~Hafstein\thanks{
email \texttt{sigurdurh@ru.is}}\\
School of Science and Engineering\\
  Reykjavik University\\
Menntavegur 1\\
IS-101 Reykjavik\\
Iceland}

\date{\today}

\maketitle

\begin{abstract}
A function $f$ from a subset of $\R^n$ to $\R$ is continuous at the origin, if and only if $\lim_{t\to 0+} f(\bs(t))=f(\bnull)$ for all continuous paths $\bs$ with $\lim_{t\to 0+} \bs(t)=\bnull$.  The continuity of $f$ can, however, be
characterized by a much smaller class of paths.  We show that the class of all paths fulfilling $\lim_{t\to 0+} \bs(t)=\bnull$, $\bs\in [\cC^\infty(]0,a[)]^n$, and $\sup_{t\in\,]0,a[}\|\bs(t)\|\|\bs'(t)\| < +\infty$ is sufficient.  Further, given any sequences $(\bx_k)_{k\in\N}$ and $(\by_k)_{k\in\N}$ in $\R^n\setminus\{\bnull\}$, such
that $\lim_{k\to+\infty}\bx_k=\bnull$, $\bx_k\cdot\by_k \ge 0$, and $\|\by_k\|=1$ for all $k\in\N$, we show that there exist a path of this class, such that $\bs(\|\bx_k\|)=\bx_k$ and $\bs'(\|\bx_k\|)=\by_k$
for an infinite number of $k\in\N$.
\end{abstract}

\section*{Background}
These results were derived by the author several years ago, in an unsuccessful attempt to work out the details of Theorem 3.2 in \cite{wilson69}, a theorem that claims the existence of smooth Lyapunov functions for
uniformly\footnote{The use of ``uniformely'' is unusual in this context and does not reefer to uniform in time, as has become widely accepted terminology, but to asymptotically stable in the sense of Lyapunov, cf.~e.g.~\cite{khalil,vidyasagar}.}, asymptotically stable, closed, invariant sets.  This claim had actually been proved earlier as Theorem 14$^*$ in \cite{massera56}, but the proof is hard to read.  Based on Massera's ideas, in a slightly more general context, a proof of this theorem is worked out in details in \cite{lin96}.

Other parts of the proof of Theorem 3.2 in \cite{wilson69} had already been clarified in \cite{nadzieja90}, but some arguments there were not comprehensible to the author. The proof in \cite{wilson69} and the 
arguments by \cite{nadzieja90} eventually turned out to be incorrect, as shown by a counterexample on page 5 in \cite{HafsteinMon}.  The attempt to work out the details was thus doomed to fail.  The lemmas in this paper were intended as auxiliary results in this attempt and are published here in the hope that they might be useful for someone else in a different context.

We denote by $\|\cdot\|$ the euclidian norm on $\R^n$, by $\N$ the set of the integers larger than zero, 
and by $\bx\cdot\by$ the scalar product of $\bx,\by\in\R^n$.
\section*{Results}
Make the following observation:
 The continuity of a function $f:\R^n\to\R$ at a point $\by$ can be characterized through
sequences.  That is, $f$ is continuous at $\by$, if and only if for every sequence $(\bx_k)_{k\in\N}$ with $\bx_k\to \by$  we have $f(\bx_k)\to f(\by)$.
Obviously ``every sequence'' can be mollified somewhat without affecting the results, e.g.~to every sequence $(\bx_k)_{k\in\N}$ fulfilling $\|\bx_k-\by\| < 1/k$ or, less obviously, for a fixed $\vartheta>0$, every sequence $(\bx_k)_{k\in\N}\subset C_\vartheta$ of every right circular cone $C_\vartheta$ with the vertex at $\by$ and $\vartheta$ as the aperture (opening angle).  The second claim can be proved almost identically to Lemma 2 below.

To show the continuity of $f$ at $\by$ one can alternatively consider continuous paths $\bs:\,]0,a,[\to\R^n$, $a>0$, such that $\lim_{t\to 0+}\bs(t)=\by$.  Then clearly $f$ is
continuous, if and only if $\lim_{t\to 0+}f(\bs(t))=f(\by)$ for all such paths.  Just as in the case of sequences the condition ``every continuous path'' can be mollified.
We mollify it to all paths $\bs\in\cC^\infty(]0,a[)$, $a>0$, such that $\sup_{t\in\,]0,a[}\|\bs(t)\|\|\bs'(t)\| < \infty$ in the next three lemmas.  We assume, without loss of
generality, that $f(\by)=0$ and $\by=\bnull$.

Before we prove the results claimed in the abstract we first prove a simple lemma.
\begin{lemma}
\label{DIFFPWLEMMA}
Let $m\in\N$, $m\ge 2$, and let $-\infty< t_1 <t_2< \ldots < t_m < +\infty$.  Assume that the function $\bp:\R\to\R^n$ is continuous and affine on every interval $[t_i,t_{i+1}]$, $i=1,2,\ldots,m-1$, i.e.~there are
vectors $\ba_i,\bb_i\in\R^n$, $i=1,2,\ldots,m-1$, such that
$$
\bp(t) := \ba_i t +\bb_ i,\ \ \ \ \text{if\ \  $t_i \leq t \leq t_{i+1}$}.
$$
Then, for every nonnegative function $\rho\in \cC^\infty(\R)$, such that $\int_\R\rho(\tau)d\tau=1$ and $\tau \mapsto \rho(t-\tau)$ has compact support in $[t_1,t_m]$,  we have
$$
\left|\frac{d}{dt}\int_\R \rho(t-\tau)\bp(\tau)d\tau\right| \leq \sum_{i=1}^{m-1}\|\ba_i\|.
$$
\end{lemma}
\Proof
By  partial integration
\begin{align*}
\frac{d}{dt} \int_\R\rho(t-\tau)\bp(\tau)d\tau &= \frac{d}{dt} \int_{t_1}^{t_{m}}\rho(t-\tau)\bp(\tau)d\tau
= \sum_{i=1}^{m-1}\int_{t_i}^{t_{i+1}}\rho'(t-\tau)(\ba_i \tau+ \bb_i)d\tau \\
&= \sum_{i=1}^{m-1} \left(  -\rho(t-\tau)(\ba_i \tau +\bb_i)\Big{|}_{\tau=t_i}^{t_{i+1}} + \int_{t_i}^{t_{i+1}}\rho(t-\tau)\ba_i d\tau\right).
\end{align*}
Because $\bp$ is continuous we have $\ba_i t_i +\bb_i = \ba_{i+1} t_i + \bb_{i+1}$ for all $i=1,2,\ldots,m-1$ and because $\tau\mapsto \rho(t-\tau)$ has compact support in $[t_1,t_m]$ we have
$\rho(t-t_1) = \rho(t-t_m)=0$.  Therefore
$$
\left|\frac{d}{dt}\int_\R \rho(t-\tau)\bp(\tau)d\tau\right| \leq \sum_{i=1}^{m-1}\|\ba_i\| \int_{t_i}^{t_{i+1}}\rho(t-\tau) d\tau \leq \sum_{i=1}^{m-1}\|\ba_i\|.
$$
\eot

\begin{lemma}
\label{LIMLEMMA}
Let $(\bx_k)_{k\in\N}$ and $(\by_k)_{k\in\N}$ be sequences in $\R^n\setminus\{\bnull\}$ with the properties that\,{\rm :}
\begin{itemize}
\item[a)] $\lim_{k\ra +\infty} \bx_k = \bnull$.
\item[b)] $\|\by_k\|=1$ for every $k\in\N$.
\item[c)] $\bx_k \cdot \by_k \geq 0$ for every $k\in\N$.
\end{itemize}
Then there exists a path $\bs:\,]0,a[\, \to \R^n$, $a>0$, with the following properties{\rm :}
\begin{itemize}
\item[i)]
$\bs \in [\cC^\infty(]0,a[)]^n$.
\item[ii)]
$\lim_{t\ra0+}\bs(t)=\bnull.$
\item[iii)]
$\sup_{t\in\,]0,a[}\|\bs(t)\|\|\bs'(t)\| < +\infty.$
\item[iv)]
The set $\{ k\in\N \,:\, \bs(\|\bx_k\|)=\bx_k\  \text{and}\ \  \bs'(\|\bx_k\|) = \by_k\}$ has an infinite number of elements.
\end{itemize}
\end{lemma}

\Proof
Set $\cX := \{\bx_k \,:\, k\in \N\}$.  Obviously $\cX$ contains an infinite number of elements.
For all $\bz \in \R^n$ with $\|\bz\|=1$ we define
$$
\cK_\bz := \left\{\bx \in \R^n \,:\, \|\bx - r\bz\| \leq \frac{r}{\sqrt{3}} \ \ \text{for some $r\geq 0$}\right\}.
$$
Then $\cK_\bz$ is a right circular cone with the vertex at the origin and the opening angle at the vertex is $\pi/3$.

There is a finite number $\cK_{\bz_1},$ $\cK_{\bz_2},\ldots,\cK_{\bz_m}$ of such sets,
such that
$$
\{\bx\in\R^n\,:\,\|\bx\| \le 1\} \subset \bigcup_{i=1}^m \cK_{\bz_i}.
$$
To see this let the set $\cS^n := \{\bx\in \R^n\,:\, \|\bx\| =1\}$ be equipped with its usual topology (see, for example, Section 8-8 in \cite{hocking}).
For every $\|\bz\|=1$ define $\cO_\bz$
to be the interior of $\cK_\bz\cap\cS^n$ in the $\cS^n$ topology.  Because $\cS^n$ is compact, a finite number of the $\cO_\bz$\,s, say
$\cO_{{\bz_1}},\cO_{{\bz_2}},\ldots,\cO_{{\bz_m}}$, suffice to cover $\cS^n$. But then, because $\cK_{\bz_i} = \bigcup_{r\geq0}r\cO_{\bz_i}$,
the sets $\cK_{\bz_1},$ $\cK_{\bz_2},\ldots,\cK_{\bz_m}$ cover $\R^n$.
Hence,  there is at least one $\bz^*$ such that $\cX\cap\cK_{\bz^*}$ contains an infinite number of elements.

Define the sets
$$
\cK_{\bz^*,k} := \left\{\bx\in\cK_{\bz^*} \,:\, \frac{1}{k+1} < \|\bx\| \leq \frac{1}{k} \right\}
$$
for all $k\in\N$.  Then at least one of the sets
$
\cX \cap \bigcup_{k=1}^{+\infty} \cK_{\bz^*,2k}$ or $\cap\bigcup_{k=1}^{+\infty} \cK_{\bz^*,2k-1}$
contains an infinite number of elements.  Without loss of generality we assume that
$\cX\cap\bigcup_{k=1}^{+\infty} \cK_{\bz^*,2k}$ is infinite.  The rest of the proof would be almost identical under the alternative assumption.

We construct sequences $(\ba_i)_{i\in\N}$ and $(\bb_i)_{i\in\N}$  in $\cK_{\bz^*}$ in the following way\,{\rm :}
For every $k\in\N$ consider the intersection $\cX\cap\cK_{\by^*,2k}$.
If it is not empty there is an $i\in \N$ such that $\bx_i$ from the sequence $(\bx_i)_{i\in\N}$ is in $\cK_{\by^*,2k}$ and
in this case we set $\ba_k := \bx_i$ and $\bb_k := \by_i$.  If the intersection is empty
we set $\ba_k$ equal to an arbitrary element of $\cK_{\by^*,2k}$ and set $\bb_k := \ba_k/\|\ba_k\|$.
By this construction there are infinitely many $i$ and $k$ in $\N$ such that $\ba_k = \bx_i$ and $\bb_k = \by_i$.

We now have everything we need to construct the claimed path $\bs$.
We start by constructing a piecewise affine path $\widetilde\bs\in\cC(]0,\|\ba_1\|])$ and then smooth it to get $\bs$.

For every $k\in\N$ we define for $m=0,1,2$ the constants $t_{k,m}$ by
$$
t_{k,m} := \|\ba_k\| - \frac{m}{3}\left(\frac{1}{2k+1} - \frac{1}{2k+2}\right).
$$
By the construction of the sequence $(\ba_k)_{k\in\N}$ we have $0<t_{k+1,0}<t_{k,2}<t_{k,1}<t_{k,0}$ for all $k\in\N$.
We define $\widetilde\bs$ for every $k\in\N$ on the interval $]t_{k+1,0},t_{k,0}]$   by
$$
\widetilde\bs(t) :=  \ba_{k+1} + (t-t_{k+1,0})\bb_{k+1}\ \ \ \text{for all $t \in\, ]t_{k+1,0},t_{k,2}]$,}
$$
$$
\widetilde\bs(t) :=\ba_{k+1} +
(t_{k,2} - t_{k+1,0})\bb_{k+1} + (t-t_{k,2})\frac{\ba_k - \ba_{k+1} +(t_{k,1}-t_{k,0})\bb_k - (t_{k,2}-t_{k+1,0})\bb_{k+1}}{t_{k,1}-t_{k,2}}
$$
for all $t\in\,]t_{k,2},t_{k,1}]$, and
$$
\widetilde\bs(t) := \ba_k + (t-t_{k,0})\bb_k\ \ \ \text{for all $t\in\,]t_{k,1},t_{k,0}]$.}
$$

Then
\begin{equation}
\label{PROP4}
\widetilde\bs(\|\ba_k\|) = \widetilde\bs(t_{k,0}) = \ba_k\ \ \text{and}\ \ \widetilde\bs'(\|\ba_k\|) = \widetilde\bs'(t_{k,0}) = \bb_k\ \ \text{for all $k\in\N$.}
\end{equation}
Further, $\widetilde\bs$ is continuous on $]0,t_{1,0}]$ and smooth, except at the
points $t=t_{k,1}$ and $t=t_{k,2}$ for all $k\in\N$.

Let $\rho\in\cC^\infty(\R)$ be a nonnegative function with $\supp(\rho)\subset [-1,1]$ and $\int_\R\rho(\tau)d\tau = 1$.
We define the smooth path $\bs:\,]0,t_{1,0}[\, \lra \R^n$ by defining, whenever
$$
\frac{t_{k+1,0} + t_{k,2}}{2} \leq t \leq \frac{t_{k,1}+t_{k,0}}{2}
$$
for some $k\in\N_{>0}$,
$$
\bs(t) := \int_{t-\frac{t_{k,0}-t_{k,1}}{4}}^{t+\frac{t_{k,0}-t_{k,1}}{4}}\rho\left(\frac{4(t-\tau)}{t_{k,0}-t_{k,1}}\right)
\frac{4\ \tilde\bs(\tau)}{t_{k,0}-t_{k,1}}d\tau = \int_{-1}^1 \rho(\tau)\ \tilde\bs\left(t -\tau\frac{t_{k,0}-t_{k,1}}{4}\right)d\tau,
$$
 and we set $\bs(t) := \tilde\bs(t)$ otherwise.  Then $\bs$ fulfills the claimed property {\it i)}.

Note that because $t_{k,0} > (2k+1)^{-1}$ and $t_{k+1,0} \leq (2k+2)^{-1}$, which implies $t_{k+1,0}-t_{k,2} \geq t_{k,0}-t_{k,1}$,
the paths $\bs$ and $\widetilde\bs$ coincide for every $t$ such that
$$
t_{k+1,0} - \frac{t_{k,0}-t_{k,1}}{4} < t < t_{k+1,0} + \frac{t_{k,0}-t_{k,1}}{4}
$$
for every $k\in\N$.  Hence the path $\bs$ fulfills the claimed property {\it iv)}.  

For every $k\in\N$ and every
$0<t<t_{k,0}$ we have the estimate
$$
\|\bs(t)\| \leq \sqrt{t_{k,0}^2 + (t_{k,0} - t_{k,1})^2},
$$
so, for every $k\in\N$ and every $0<t<(2k)^{-1}$ we have the crude estimate
\begin{equation}
\label{PROP2}
\|\bs(t)\| \leq \frac{1}{k}
\end{equation}
and the claimed property {\it ii)} is fulfilled as well.  We come to the claimed property {\it iii)}.
By Lemma \ref{DIFFPWLEMMA} we have for every $k\in\N$ and every $t_{k+1,0} < t \leq t_{k,0}$ the estimate
\begin{align*}
\|\bs'(t)\| &\leq \|\bb_k\| + \|\bb_{k+1}\| +
\left\|\frac{\ba_k - \ba_{k+1} +(t_{k,1}-t_{k,0})\bb_k - (t_{k,2}-t_{k+1,0})\bb_{k+1}}{t_{k,1}-t_{k,2}}\right\| \\
&\leq 2 +3(2k+1)(2k+2)\left(\frac{k+1}{k(2k+2)}  + \frac{1}{k(2k+2)} - \frac{1}{3(2k+1)(2k+2)}\right)\\
&\leq 28k,
\end{align*}
where we used
$$
\|\ba_k - \ba_{k+1}\| = \sqrt{\|\ba_k\|^2 + \|\ba_{k+1}\|^2 - 2\ba_k\cdot\ba_{k+1}},
$$
which has a maximum with $\|\ba_k\| = (2k)^{-1}$, $\|\ba_{k+1}\| = (2k+2)^{-1}$, and $\ba_k\cdot\ba_{k-1} = \|\ba_k\|\|\ba_{k+1}\|/2$
for $k\geq 2$
(recall that the opening angle at the vertex of $\cK_{\bz^*}$ is $\pi/3$).
But then
$$
\|\bs(t)\|\|\bs'(t)\| \leq 28
$$
for all $t\in\,]0,t_{1,0}[$ and we have finished the proof.
\eot

\begin{lemma}
\label{LIMTHEO}
Let $f:\cN \to \R$, where $\cN\subset \R^n$ is a neighbourhood of the origin.
Assume $f(\bnull)=0$.  Then $f$ is continuous at the origin, if and only if 
for every path $\bs\in[\cC^\infty(]0,a[)]^n$, $a>0$, such that
\begin{equation}
\label{VFC}
\lim_{t\ra0+}\bs(t)=\bnull\ \ \ \text{and}\ \ \ \sup_{t\in\,]0,a[}\|\bs(t)\|\|\bs'(t)\| < +\infty,
\end{equation}
we have
$$
\lim_{t\ra 0+} f(\bs(t))= 0.
$$
\end{lemma}
\Proof
The ``only if'' part is obvious.  We prove the ``if'' part by showing that if $f$ is not continuous at the origin, then there is a path
$\bs\in[\cC^\infty(]0,a[)]^n$, $a>0$, that fulfills the properties (\ref{VFC}), but for which $\limsup_{t\ra 0+} |f(\bs(t))| > 0$.

Assume that $f$ is not continuous at the origin.  Then there is
an $\varepsilon>0$ and a sequence $\bx_k$, $k\in\N$, such that $\lim_{k\ra+\infty}\bx_k = \bnull$ but $|f(\bx_k)| \geq \varepsilon$ for all $k\in\N$.
Set e.g.~$\by_k:=\bx_k/\|\bx_k\|$ for all $k\in\N$.
By Lemma \ref{LIMLEMMA} there is a path $\bs\in[\cC^\infty(]0,a[)]^n$, $a>0$, with the properties that
$\lim_{t\ra 0+} \bs(t) = \bnull$, $\sup_{t\in\,]0,a[}\|\bs(t)\|\|\bs'(t)\| < +\infty$,
and $\bs(\|\bx_k\|) = \bx_k$ for an infinite number of $k\in\N$.
But then $|f(\bs(\|\bx_k\|))| \geq \varepsilon$ for an infinite number of $k\in\N$, which implies
$\limsup_{t\ra 0+} |f(\bs(t))| \geq \varepsilon$ and we have finished the proof.
\eot

\end{document}